\documentclass[11pt]{amsart}
\usepackage{graphicx} % Required for inserting images

\usepackage{amsfonts}
\usepackage{amsmath}
\usepackage{amsthm}
\usepackage{latexsym}
\usepackage{amssymb}

\usepackage{amscd}
\usepackage{xcolor}
\usepackage{hyperref}
\usepackage{alphabeta}
\usepackage{fancyhdr}
\usepackage{mathtools}
\usepackage{fancyvrb}

\usepackage{imakeidx}
\usepackage{calrsfs}
\usepackage{soul}
\usepackage{multirow}
\usepackage{todonotes}
\usepackage{amsbsy}
\usepackage[all]{xy}
\usepackage{changes}
\usepackage{enumerate}
\usepackage{comment}
\usepackage{mathtools}
\usepackage{float}
\usepackage{color, colortbl}
\usepackage{changes}
\usepackage{makeidx}

\theoremstyle{definition}

\theoremstyle{remark}

\numberwithin{equation}{section}

\newcommand{\abs}[1]{\lvert#1\rvert}

\newcommand{\norm}[1]{\lVert#1\rVert}

\def\D{{\mathbb D}}  \def\T{{\mathbb T}}
\def\C{{\mathbb C}}

\def\bD{{\mathbb D}}   \def\cB{{\mathcal B}}

\def\cH{{\mathcal H}}

\begin{document}

\title{ On the Hilbert matrix operator: a brief survey}

\author{C. Bellavita }
\email{carlo.bellavita@gmail.com}
\address{Department of Mathematics, Aristotle University of Thessaloniki, 54124, Thessaloniki, Greece.}
\address{Current address: Departamento de Matem\'atica i Inform\'atica, Universitat de Barcelona, Gran Via 585, 08007 Barcelona, Spain.
}

\author{V. Daskalogiannis}
\email{vdaskalo@math.auth.gr}
\address{Department of Mathematics, Aristotle University of Thessaloniki, 54124, Thessaloniki, Greece}
\address{Division of Science and Technology, American College of Thessaloniki, 17 V. Sevenidi St., 55535, Pylea, Greece.}

\author{G. Stylogiannis}
\email{stylog@math.auth.gr}
\address{Department of Mathematics, Aristotle University of Thessaloniki, 54124, Thessaloniki, Greece.}


\thanks{This review article was supported by the
Hellenic Foundation for Research and Innovation (H.F.R.I.) under the '2nd Call for H.F.R.I. Research Projects to support Faculty Members \& Researchers' (Project Number: 4662).}

\begin{abstract}
%The aim of this article is to unfold the recent history behind the study of the Hilbert matrix acting as an operator on spaces of analytic functions, as well as on sequence spaces. We present the most recent developments on the topic, hoping that this survey will be helpful for researchers willing to dive into the captivating theory of operator matrices.
This article aims to explore the most recent developments in the study of the Hilbert matrix, acting as an operator on spaces of analytic functions and sequence spaces. We present the latest advances in this area, aiming to provide a concise overview for researchers interested in delving into the captivating theory of operator matrices.
\end{abstract}

\maketitle

\section{Introduction}
\noindent The infinite matrix
\[\label{eq:matrix}
\mathbb{H} \,=\,
\left[
\begin{array}{ccccc}
1                 & \frac{1}{2}    &\frac{1}{3}      &\frac{1}{4}       & \cdots
\\[0.1in]
\frac{1}{2}       & \frac{1}{3}    &  \frac{1}{4}    &  \frac{1}{5}     &  \cdots
\\[0.1in]
\frac{1}{3}       &\frac{1}{4}     & \frac{1}{5}     & \frac{1}{6}      & \cdots

\\[0.1in]
\vdots            &\vdots          & \vdots          & \vdots           & \ddots 
\end{array}
\right]\,=\,\left(\frac{1}{n+k+1}\right)_{n,k\,=\,0,1,2...}
\]
is the well known Hilbert matrix. It has been an object of study for more than a century and it has several applications in different fields of science: see for example \cite{Collar} for aerodynamics, \cite{Magnus} for diffraction of electromagnetic waves and \cite{Lukacs} for statistics.

This matrix was first introduced by David Hilbert in 1894 \cite{Hilbert1900}, in connection to {\it Legendre} polynomials and the {\it least squares} approximation theory. 
It turns out that Hilbert's matrix is strongly ill-conditioned and that even its finite $n\times n$ sections, denoted by $\mathbb{H}_n$, “almost fail" to be invertible for large values of $n$. Indeed, Hilbert noticed that the determinant of each finite section is the reciprocal of an integer. In particular the inverse of $\mathbb{H}_n$ has very large integer entries, making the computation of its eigenvalues a very sensitive problem. 
 As a measure of comparison, we note that $\det{(\mathbb{H}_9)}\simeq 10^{-42}$. Hilbert's matrix is symmetric, {\it positive definite}, meaning that all its eigenvalues are positive, and {\it totally positive}, since every submatrix has a positive determinant.

The Hilbert matrix is also the canonical example of a {\it Hankel} matrix: we recall that an infinite matrix $A=(a_{i,j})$ is a Hankel matrix if its entries depend only on the sum of the indices, that is if $a_{i,j}=a_{i+j}$. For more information on Hankel operators, the interested reader is referred to V. Peller's monograph \cite{Peller2003}.
 The Hilbert matrix is, without a doubt, one of the most important examples of such matrices, related to a number of different areas of mathematics such as Numerical analysis, Number theory and Linear algebra. Furthermore, it also plays a prominent role in the areas of Operator theory and Complex analysis, where we focus in the sequel.
%%%%%%%%%%%%%%%%%%%%%%%%%%%%%%%%%%%%%%%%%%%%%%%%%%%%%%%%%%%%%%
\section{The Hilbert matrix acting on sequence spaces}
\noindent At the turn of the 20th century, during his lectures on integral equations, Hilbert introduces a remarkable bi-linear form and states his famous “double-series theorem", i.e.
\begin{equation}\label{Hilbert inequality 1}
    \sum_{m=1}^{\infty} \sum_{n=1}^{\infty}\dfrac{a_m b_n}{m+n}\leq
C \left(\sum_{m=1}^{\infty}{a_m}^2\right)^{\frac{1}{2}} \left(\sum_{n=1}^{\infty} {b_n}^2\right)^{\frac{1}{2}}\,.
\end{equation}
As G. H. Hardy, J. E. Littlewood and G. P\'{o}lya describe in \cite{Hardy1934}, Hilbert's proof was published by H. Weyl in his dissertation \cite{Weyl1908}. The determination of the exact constant $C=\pi$ is due to I. Schur in 1911 \cite{Schur1911}. Hardy and M. Riesz generalized \eqref{Hilbert inequality 1} to other values of $p\neq 2$ \cite{Hardy1925} and, furthermore, Hardy proved a sharper, non-homogeneous form \cite[Theorem 323]{Hardy1934}, now known as {\it Hilbert's inequality:}
\begin{equation}\label{Hilbert inequality}
 \sum\limits_{m=0}^{\infty} \sum\limits_{n=0}^{\infty}\dfrac{a_m b_n}{m+n+1}\leq
\dfrac{\pi}{\sin{\left(\frac{\pi}{p}\right)}} \left(\sum\limits_{m=0}^{\infty}{a_m}^p\right)^{\frac{1}{p}} \left(\sum\limits_{n=0}^{\infty} {b_n}^q\right)^{\frac{1}{q}} \,,  
\end{equation}
where $\{a_m\},\,\{b_n\}$ are sequences of positive terms, with $p,q>1$ such that $\frac{1}{p}+\frac{1}{q}=1$, and the constant ${\pi}/{\sin{\frac{\pi}{p}}}$ is the best possible.

Considering any sequence as a column vector and applying formal matrix multiplication, we see that the Hilbert matrix induces a sequence transformation in the following way:
\[
\{a_n\}_n \longmapsto \left\{\sum_{k=0}^\infty\dfrac{a_k}{n+k+1}\right\}_n\, ,
\]
for those $\{a_n\}$ for which the right-hand sums are well defined.
Hilbert's inequality readily implies that the Hilbert matrix induces a bounded linear operator, denoted by $\cH$, which acts on the classical space of $p$-summable sequences $\ell^p$.
By the duality of the spaces $\ell^p$, inequality \eqref{Hilbert inequality} is equivalent to
\begin{equation}\label{dual}
\left(\sum_{n} \left|\sum_k \dfrac{a_k}{n+k+1} \right|^p \,\right)^{\frac{1}{p}}\, \leq\, \dfrac{\pi}{\sin{\frac{\pi}{p}}}\; \norm{\{a_n\}}_{\ell^p}\,.
\end{equation}
Consequently $\cH: \ell^p \to \ell^p,\;\;1<p<\infty$, is bounded thus $$\norm{\cH}_{\ell^p \to \ell^p}=\frac{\pi}{\sin{\frac{\pi}{p}}}\,.$$

\noindent
Quite recently, the authors in \cite{Daskalogiannis2023} showed that $\cH$ induces a bounded operator on the spaces $\ell_{p-2}^p$ when $ 1<p<\infty$, with the same operator norm ${\pi}/{\sin{\frac{\pi}{p}}}$ as in the unweighted case. We recall that
$$
\norm{\{a_n\}}^p_{\ell_{p-2}^p}=\sum_{n=0}^\infty (n+1)^{p-2} \abs{a_n}^p <\infty\,.
$$
It is worth mentioning that the method employed in \cite{Daskalogiannis2023} can be adapted to offer an alternative proof of \eqref{Hilbert inequality}. To conclude this section, we also note that the action of the Hilbert matrix between various sequence spaces, including the Césaro matrix domain, Copson sequence spaces, and Gamma sequence spaces was investigated by H. Roopaei in \cite{Roopaei2020}.

%Closing this section, we mention that H. Roopaei \cite{Roopaei2020} studied  the action of the Hilbert matrix between various sequence spaces, such as the C\'esaro matrix domain, Copson sequence spaces and Gamma sequence spaces.
%%%%%%%%%%%%%%%%%%%%%%%%%%%%%%%%%%%%%%%%%%%%%%%%%%%%%%%%%%%%%%
\section{The Hilbert matrix as an operator on spaces of analytic functions}
\noindent In 2000, E. Diamantopoulos and A. G. Siskakis \cite{Diamantopoulos2000} initiated the study of the Hilbert matrix acting as an operator on spaces of analytic functions. They considered its action on the sequence of the corresponding Taylor coefficients, inducing, in this way, an operator acting on the Hardy spaces of the unit disc.

Let $\D$ denote the unit disc of the complex plane, let $\T$ be its boundary and denote by $H(\D)$ the space of all holomorphic functions in $\D$. The Hardy space $H^p$ consists of all $f\in H(\D)$, such that
\[
\norm{f}_{H^p}= \sup_{0<r<1} M_p(r, f) <\infty\,,
\]
where $M^p_p(r,f) = \dfrac{1}{2\pi} \int_0^{2\pi} \abs{f(re^{it})}^p\,dt$ are the integral means of $f$, along circles of radius $r$. The space $H^2$ coincides with $\ell^2$ and hence it is a Hilbert space, while for $p<q$ we have that $H^q \subsetneq H^p$.

For a holomorphic function $f(z)=\sum_n a_n\,z^n \in H^1$, define 
\begin{equation}\label{Hf}
\cH(f)(z):=\sum\limits_{n=0}^{\infty} \sum\limits_{k=0}^{\infty}\dfrac{a_k}{n+k+1}\,z^n\,.
\end{equation}
Hardy's inequality \cite[p. 48]{Duren1970}
\[
\sum_n \dfrac{\abs{a_n}}{n+1}\leq \pi \norm{f}_{H^1}
\]
implies that the power series defining $\cH(f)$ has bounded coefficients; therefore, its radius of convergence is greater than or equal to 1, and $\cH(f)$ defines an analytic function of the unit disc for every $f\in H^1$.
Notice that the series in \eqref{Hf} is not always well defined when $f$ is an arbitrary analytic function. Consider for example the function $f(z)=\frac{1}{1-z}=\sum_n z^n$, and observe that, in that case, the inner sum in \eqref{Hf} diverges. This means that one needs to restrict the domain of the operator on spaces for which the right hand side of \eqref{Hf} makes sense.

Diamatopoulos and Siskakis \cite{Diamantopoulos2000} studied the boundedness of the operator
\[
\cH: H^p \to H^p,\;\;1\leq p \leq \infty\,.
\]
Due to the nature of these spaces, they obtained an integral representation for $\cH$, that is
\[
\begin{split}
\cH(f)(z)&=\sum_{n=0}^{\infty} \sum_{k=0}^{\infty}\dfrac{a_k}{n+k+1}\,z^n
=\sum_{n=0}^{\infty}\sum_{k=0}^{\infty} a_k\,\int_0^1 t^{n+k} \,dt\, z^n
\\
&=\sum_{n=0}^\infty\int_{0}^{1}f(t)\, t^n\, z^n\,dt= \int_{0}^{1}f(t)\sum_{n=0}^\infty (t z)^n\,dt\,,
\end{split}
\]
hence
\begin{equation}\label{Integral form 1}
\cH(f)(z)=\,\int_{0}^{1}f(t)\dfrac{1}{1-tz}\,dt\,.  
\end{equation}
The convergence of all sums and integrals involved is guaranteed by Hardy's inequality and by the Fej\'er-Riesz inequality \cite[Theorem 3.13]{Duren1970} for functions in $H^p$,
\[
\int_{-1}^1 \abs{f(t)}^p\,dt \leq \frac{1}{2}\norm{f}^p_{H^p}\,,
\]
which allows the interchange between sums and integrals. In \cite{Diamantopoulos2000},  it was shown that $\cH$ is bounded on $H^p,\;1<p<\infty$, but not bounded in $H^1$ and $H^\infty$. For $H^\infty$, it suffices to consider the constant function $f(z)=1$. Then, $\cH(1)=\frac{1}{z}\log\frac{1}{1-z}$, which clearly is not bounded. For the case $p=1$, consider  the function $f_\epsilon(z)=(1-z)^{-1}(\frac{1}{z}\log\frac{1}{1-z})^{-1-\epsilon}$, which is in $H^1$ for each $\epsilon > 0$ \cite[p.13]{Duren1970}. Its image under $\cH$ is given by
\[
\cH(f_\epsilon)(z)=\sum_{n=0}^\infty \int_0^1 f_\epsilon(t) t^n\,dt\,z^n
\]
which is analytic on $\mathbb D$.
Assuming that it is an $H^1$ function, Hardy's inequality implies that the quantity
\[
\int_0^1 \dfrac{1}{(1-t)(\frac{1}{t}\log\frac{1}{1-t})^{\epsilon}}\,dt
\]
is finite. This leads to a contradiction when $\epsilon \leq 1$.

For $p>1$, one can change the path of integration in \eqref{Integral form 1}, to be the arc joining $0$ and $1$, given by
\[
\gamma_z(t)=\dfrac{t}{1-(1-t)z},\;\;0\leq t \leq 1\,.
\]
Then,
\[
\begin{split}
\cH(f)(z) &= \int_0^1 \dfrac{1}{1-(1-t)z}\,f\left(\dfrac{t}{1-(1-t)z} \right)\,dt
\\
&=\int_0^1 w_t(z)\,f(\phi_t(z))\,dt\,,
\end{split}
\]
and $\cH(f)$ is an {\it average} of the weighted composition operators
\begin{equation}\label{Tt}
T_t(f)(z)=w_t(z) f(\phi_t(z))\,.
\end{equation}
To this point, by an application of the generalized Minkowski's inequality, one can obtain an upper bound for the norm of $\cH$, in terms of the norm of $T_t$. In specific
\begin{equation}\label{Minkowski}
    \norm{\cH(f)}_{H^p} \leq \int_0^1 \norm{T_t(f)}_{H^p}\,dt\,.
\end{equation}
An involved, technical calculation (see \cite[Section 3]{Diamantopoulos2000}), gives that for $2\leq p<\infty$
\[
\norm{T_t(f)}_{H^p}\leq \dfrac{t^{\frac{1}{p}-1}}{(1-t)^{\frac{1}{p}}}\,\norm{f}_{H^p}\,,
\]
and a less sharp upper bound for the case $1<p<2$. Using \eqref{Minkowski}, we get \cite[Theorem 1.1]{Diamantopoulos2000} that for $p\geq 2$
\[
\begin{split}
 \norm{\cH(f)}_{H^p} &\leq    \int_0^1 \dfrac{t^{\frac{1}{p}-1}}{(1-t)^{\frac{1}{p}}}\,dt\,\norm{f}_{H^p}
 \\
 &=B(\frac{1}{p},1-\frac{1}{p})\,\norm{f}_{H^p}
 \\
 &=\dfrac{\pi}{\sin{\frac{\pi}{p}}}\,\norm{f}_{H^p}\,,
\end{split}
\]
where $B(x,y)$ is the classical Beta function (see \cite{Lebedev1965}).

A few years later, in 2008, M. Dostani\'c, M. Jevti\'c, and D. Vukoti\'c \cite{Dostanic2008}, accomplished remarkable progress, clearing the scene for the Hardy spaces. We should mention here, that in the meantime, a  decisive step forward was made by B. Hollenbeck and I. Verbitsky \cite{HOLLENBECK2000}, who proved a long-standing conjecture (see \cite{Gokhberg1968}) of I. Gohberg and N. Krupnik (1968), regarding the exact value of the norm of the Riesz projection $P_+$. Remember that $P_+$ acts on a complex-valued $L^p$ function $f$ on the unit circle
\[
\begin{split}
&\quad P_+: L^p(\T)\to H^p 
\\
&\sum_{n=-\infty}^\infty a_n e^{in\theta} \longmapsto \sum_{n=0}^\infty a_n e^{in\theta}\,,
\end{split}
\]
where $a_n$ are the Fourier coefficients of $f$.
In specific, they proved that $\norm{P_+}_p={1}/{\sin\frac{\pi}{p}},\;\;1<p<\infty$.

Dostani\'c et al., considering $\cH$ in the more general setting of Hankel operators, exploited the relation
\begin{equation}\label{Nehari}
\cH=P_+ \circ M_\varphi \circ J\,,
\end{equation}
where $J$ is the isometric conjugation operator (or {\it Flip operator}), and $M_\varphi$ is the multiplication operator by the function $\varphi(t)=i(\pi-t)e^{-it} \in L^\infty (\T)$. Consequently
\[
\norm{\cH(f)}_{H^p}\leq \norm{P_+}_p\, \norm{\varphi}_\infty \, \norm{Jf}_{L^p(\T)}\,,
\]
or equivalently
\[
\norm{\cH(f)}_{H^p}\leq \dfrac{\norm{\varphi}_\infty}{\sin\frac{\pi}{p}}\;\norm{f}_{H^p}=\,\frac{\pi}{\sin\frac{\pi}{p}}\;\norm{f}_{H^p},\;\;1<p<\infty\,.
\]
Furthermore, they were able to obtain a sharp lower bound, using the test functions $f_\gamma(z)=(1-z)^{-\frac{\gamma}{p}},\;0<\gamma<1$. Some demanding calculations give
\[
\cH(f_\gamma)(z)=g(z)-R(z)\,,
\]
where $R$ has a uniformly bounded $L^p(\T)$ norm and the function $g$ satisfies $\norm{g}_{L^p(\T)}=\dfrac{\pi}{\sin\frac{\pi\gamma}{p}}\,\norm{f_\gamma}_{H^p}$. This implies that
\[
\norm{\cH}_{H^p\to H^p} \norm{f_\gamma}_{H^p}\geq \norm{\cH(f_\gamma)}_{H^p}\geq \,\left|\norm{g}_{L^p(\T)}-\norm{R}_{L^p(\T)}\right|\,
\]
or equivalently
\[
\norm{\cH}_{H^p\to H^p} \geq \left| \dfrac{\pi}{\sin\frac{\pi\gamma}{p}}- \dfrac{\norm{R}_{L^p(\T)}}{\norm{f_\gamma}_{H^p}} \right|\,.
\]
Letting $\gamma \to 1^-$, and taking into account that $\norm{f_\gamma}_{H^p}\to \infty$, the result follows (\cite[Theorem 2]{Dostanic2008}). Hence, for $1<p<\infty$, we have that (\cite[Corollary 1]{Dostanic2008})
\[
\norm{\cH}_{H^p\to H^p}=\dfrac{\pi}{\sin\frac{\pi}{p}}\,.
\]
In 2012, B. Lanucha, M. Nowak and M. Pavlovi\'c \cite{Lanucha2012} mainly focused on the action of the Hilbert matrix on the “extreme" cases $H^1$ and $H^\infty$. They proved, among other results, that if $f\in H^1$, then $\cH(f)$ extends to a continuous function on $\overline{\D}\setminus\{1\}$ and $\cH: H^1 \to H^p,\;0<p<1$ is bounded. Furthermore, they prove \cite[Theorems 2.3, 2.4]{Lanucha2012} that if $f\in H^1$ is such that
\[
\int_{-\pi}^\pi \abs{f(e^{it})} \log \frac{\pi}{t}\,dt <\infty\,,
\]
then $\cH(f) \in H^1$, and in the special case where $f$ has positive Taylor's coefficients, then $\cH(f) \in H^1$ if and only if 
\[
\sum_{n=0}^\infty \dfrac{a_n \log(n+2)}{n+1}<\infty \,.
\]
They also notice that since the space BMOA
%, consisting of the functions in $H^1$ whose {\it boundary function} is of bounded mean oscillation on $\T$, 
is the Riesz (Szeg\"o) projection of $L^\infty (\T)$ \cite[Theorem 9.21]{Zhu1990}, then exploiting \eqref{Nehari}, one has that $\cH : H^\infty \to BMOA$ is bounded.
 Additionally, in \cite[Theorem 1.1]{BellavitaS2024}, the authors proved that
\[
\|\mathcal{H}\|_{H^\infty\to BMOA}=1+\frac{\pi}{\sqrt{2}}\, .
\]
 Very recently, Y. Guo and P. Tang \cite{guo2024} showed that the range of $\cH$ acting on $H^\infty$ is contained in the Zygmund-type space $\mathcal{Z}_1$, consisting of functions $f$ for which
 \[
 \sup_{0<r<1}(1-r^2)M_1(r, f'')<\infty\,.
 \] They also proved that $\norm{\cH}_{H^\infty \to \cB}=3$, where $\cB$ is the classical Bloch space.

Furthermore, the authors in \cite{Lanucha2012}, also considered the action of $\cH$ on the weighted Hardy spaces $H^p_\alpha$, defined by
\[
H^p_\alpha = \left\{f\in H(\D):\; M_p(r,f)=O\left( (1-r)^{-\alpha} \right) \right\}\,,
\]
and proved that $\cH:\,H^p_\alpha \to H^p_\alpha\,$ is bounded, if and only if $\alpha+\frac{1}{p}<1$ \cite[Theorems 3.1, 3.8]{Lanucha2012}.

The study of the Hilbert matrix, was subsequently extended 
to include the Bergman spaces of the disc. For $0< p <\infty$ and $\alpha > -1$, the space $A^{p}_\alpha$ consists of all $f\in H(\D)$ for which
$$
\|f\|_{A^{p}_\alpha}:=\left( (\alpha +1)\int_{\mathbb{D}} |f(z)|^p(1-|z|^2)^\alpha dA(z)\right)^{1/p}<\infty\, ,
$$
where $dA(z)=\frac{dx dy}{\pi}$ is the normalized Lebesgue area measure. For $\alpha=0$, we write $A^p_0=A^p=L^p(\D) \cap H(\D)$. For more information on Bergman spaces, see \cite{Duren2004}.

In 2004, Diamantopoulos \cite{Diamantopoulos2004} studied the action of $\cH$ on $A^p$, proving analogous results, as in the Hardy spaces. Let us comment here that by \cite[Lemma 4.1]{Nowak2010} 
\[
\sum_{n=0}^\infty \frac{\abs{a_n}}{n+1}<\infty
\]
for all $f\in A^p,\;2<p<\infty$, hence $\cH(f)$ defines an analytic function of $\D$. However, $\cH$ is not well defined on $A^2$. It suffices to observe that if $f(z)=\sum_n\frac{1}{\log(n+1)} z^n \in A^2$, then the series defining $\cH(f)(0)$ is divergent. Nonetheless, for $p>2$, $\cH$ admits the same integral representation, as an average of the weighted composition operators $T_t$.

Diamantopoulos showed that $\cH: A^p \to A^p,\;2<p<\infty$, is bounded. Furthermore, if $p\geq 4$, then
\[
\norm{\cH(f)}_{A^p}\leq \dfrac{\pi}{\sin{\frac{2\pi}{p}}}\,\norm{f}_{A^p}\,,
\]
while for $2<p<4$ he obtained a less precise estimate of the norm \cite[Theorem 1]{Diamantopoulos2004}. Dostani\'c et al. estimated a lower bound for the norm \cite[Theorem 4]{Dostanic2008} for all $p>2$, using the same test functions $f_\gamma$, which belong to $A^p$ if and only if $\gamma<2$. In specific,
\[
\norm{\cH}_{A^p \to A^p}\geq \dfrac{\pi}{\sin{\frac{2\pi}{p}}}\,,\;2<p<\infty\,
\]
and consequently
\[
\norm{\cH}_{A^p \to A^p} = \dfrac{\pi}{\sin{\frac{2\pi}{p}}},\;4\leq p<\infty\,.
\]
The case $2<p<4$ remained open until 2017, when V. Bozin and B. Karapetrovi\'c \cite{Bozin2018} estimated the exact value of the norm $\pi/\sin{\frac{2\pi}{p}}$, for all $p>2$. They exploited the monotonicity of integral means, and obtained some improved estimates for Beta functions. Their estimates involved Sturm's theorem on the zeros of sequences of polynomials and several properties of special functions, thus leading to a very technical proof. In \cite{Lindstrom2019}, M. Lindstr\"{o}m, S. Miihkinen, and N. Wikman simplified this proof, avoiding the use of Sturm's theorem and applying some new, sharper estimates for the Beta function.

In the setting of the weighted Bergman spaces $A^p_\alpha$, things are way more complicated, and the story is far from over at this point. We will present the existing results briefly. 

The study of $\cH$ on $A^p_\alpha$ originates in \cite{Galanopoulos2014}. In 2016, Jevti\'c and Karapetrovi\'c proved that $\cH$ is bounded in $A^p_\alpha$, if and only if $1<\alpha +2 < p$ \cite{JEVTIC2017}. Karapetrovi\'c estimated the norm from below \cite[Theorem 1.1]{KARAPETROVIC2017} using properties of Gauss' hypergeometric function and techniques similar to the ones used in the classical Bergman space case, and obtained
\[
\|\mathcal{H}\|_{A^p_\alpha \to A^p_\alpha }\geq \frac{\pi}{\sin\left(\frac{(\alpha+2)\pi}{p}\right)},\;\;1<\alpha+2<p\,. 
\]
Moreover, he concluded \cite[Corollary 1.3]{KARAPETROVIC2017} that if $\alpha>0$ and $p\geq 2(\alpha+2)$, then 
$$
\|\mathcal{H}\|_{A^p_\alpha \to A^p_\alpha }= \frac{\pi}{\sin\left(\frac{(\alpha+2)\pi}{p}\right)}\, ,
$$
and conjectured that this is the case for all $p>\alpha+2$.

From this point on, subsequent studies concentrated on estimating a sharp upper bound for the “missing" range $\,\alpha+2\,<\,p\,<\,2(\alpha+2)\,$, with
 $\alpha>0$. Lindstrom et al. \cite{Lindstrom2021} prove that 
the conjecture holds for positive $\alpha$ and 
\begin{equation}\label{Bound 1}
\alpha + 2 + \sqrt{\alpha^2+\frac{7}{2}\alpha +3}\,\leq\, p\,<\,2(\alpha+2)\,.
\end{equation}
Karapetrovi\'c in \cite{Karapetrovic2020}, improves this range to
\begin{equation}\label{Bound 2}
 \alpha+2+\sqrt{(\alpha+2)^2-\left(\sqrt{2}-1/2\right)(\alpha+2)} \,\leq\,  p \,<\, 2(\alpha+2)  \,,
\end{equation}
and moreover he obtains a rough upper bound for the norm of $\cH$ on {\it negatively indexed} weighted Bergman spaces. This upper bound was slightly improved in \cite{Bralovic2022}, for the case $-\frac{1}{2}<\alpha <0$.

In 2022, Lindstr\"om, Miihkinen and D. Norrbo \cite{LINDSTROM2022} proved , for $\alpha>0$ and all $p>2+\alpha$, the following: If $\mathcal{H}$ does not attain
its norm, i.e. if 
$\| \mathcal{H}f\|_{A^p_\alpha} <\|\mathcal{H}\|_{A^p_\alpha \to A^p_\alpha } \,\norm{f}_{A^p_\alpha}$
for every $f \in A^p_\alpha$, then
\[
\|\mathcal{H}\|_{A^p_\alpha \to A^p_\alpha }= \frac{\pi}{\sin\left(\frac{(\alpha+2)\pi}{p}\right)}\, .
\]
Recently, in 2023, D. Dmitrov\'ic and Karapetrovi\'c in \cite{Dmitrovic2023} proved the conjecture for $\alpha> 0$ and
\begin{equation}\label{Bound 3}
 \frac{3\alpha}{4}+2+\sqrt{\left(\frac{3\alpha}{4}+2\right)^2-\frac{\alpha + 2}{2} } \,\leq\,  p \,<\, 2(\alpha+2)  \,.
\end{equation}
The comparison between the lower bounds in the range of $p$, in \eqref{Bound 1}, \eqref{Bound 2}, and \eqref{Bound 3}, gives an improvement at each next step.

The latest results concerning the norm of the Hilbert matrix operator in $A^p_\alpha$ were established by J. Dai in \cite{Dai2024}.
The author verified the conjecture about the norm for the cases $\alpha=1$ and $0< \alpha \leq 1/47$. Moreover, he also considered negative values of alpha and showed that the expected value of the norm holds when 
\[
-1 < \alpha < 0\, \,  \text{ and }\, \,  p \geq  2(\alpha + 2) \,. 
\]
The Hilbert matrix, also induces a bounded operator on the Korenblum growth spaces $\mathcal{A}^{-\alpha}$, for $0<\alpha<1$. In \cite{Lindstrom2019}, it was shown that the norm is equal to $\pi/\sin{(\pi \alpha)}$ when $0<\alpha<2/3$, while in \cite{Dai2022} it was proved that this estimate fails when $\alpha$ gets close to $1$, and the norm is strictly greater than this quantity. 
 In \cite{hu2024}, H. Hu and S. Ye studied the action of the Hilbert matrix operator on the logarithmically weighted Korenblum space $H^\infty_{\alpha,\log}$ and they estimated its norm.

Concluding this section, the authors would like to mention \cite[Chapter 14]{Jevtic2016}, where the action of the Hilbert matrix on Hardy and Bergman spaces is also discussed.
%%%%%%%%%%%%%%%%%%%%%%%%%%%%%%%%%%%%%%%%%%%%%%%%%%%%%%%%%%%%%
\section{Other integral representations}
\noindent
In most of the aforementioned results the authors choose to use the integral representation of $\cH$ in terms of the weighted composition operators $T_t$, when trying to estimate the norm of $\mathcal{H}$ from above, and the reason is obvious if we consider \eqref{Minkowski}. But, in fact, one can choose many different paths of integration in
\[
\cH(f)(z)=\,\int_{0}^{1}f(t)\dfrac{1}{1-tz}\,dt\,. 
\]
Let $h:\D\to \C$ be a univalent and starlike function, with $h(0)=0$. One can also consider a spiral-like function $h$, but then some mild modifications in what follows are needed. Let $\varphi_t(z)=h^{-1}(th(z)),\;\;0\leq t\leq 1$, which is a self-map of $\D$ for each $t$. Set $\psi_t(z)=\varphi_t(z)/z$, which is differentiable as a function of $t$. Since $\varphi_t(0)=0$, by Schwarz's lemma, the function 
$\psi_t$ maps $\D$ into itself and the path $\beta(t)=\beta_z(t)=\psi_t(z),\;0\leq t\leq 1$ is inside the disc, and joins $0$ to $1$. Evaluating the integral representation \eqref{Integral form 1} of $\cH$, along the path $\beta$, gives
\[
%\begin{split}
\mathcal{H}(f)(z)
%\dfrac{h(z)}{z}\int_0^1 w(\phi_t(z)) f(\psi_t(z))\,dt
%\
=\dfrac{h(z)}{z}\int_0^1 w(z\psi_t(z)) f(\psi_t(z))\,dt\, ,
%\end{split}
\]
where the weight function $w$ is given by
\[
w(z)=\dfrac{1}{(1-z)h'(z)}\,.
\]
Thus, each starlike (spiral-like) $h$ induces a representation for the operator $\cH$. The choice $h(z)=\frac{z}{1-z}$ gives the classical representation. As another example, for $h(z)=\log\frac{1}{1-z}$, we obtain $\varphi_t(z)=1-(1-z)^t$, $w(z)=1$ and consequently we get
\[
\cH(f)(z)=\dfrac{1}{z}\log\dfrac{1}{1-z} \int_0^1 f\left(\dfrac{1-(1-z)^t}{z}\right)\,dt\,.
\]
An appropriate choice of the function $h$ could provide insightful estimates for the norms of the weighted composition operators
\[
S_t(f)(z)=w(z\psi_t(z)) f(\psi_t(z))\,,
\]
and, in view of \eqref{Minkowski}, also information about the norm of $\cH$.

%%%%%%%%%%%%%%%%%%%%%%%%%%%%%%%%%%%%%%%%%%%%%%%%%%%%%%%%%%%%%%
\section{Generalized Hilbert operators}
\noindent
Over the past years, different generalizations of the Hilbert matrix operator emerged in the literature. Here we choose to highlight some of them.

Let $\mu$ be a finite positive Borel measure on $[0,1)$; for the sake of clarity, we write $\mu([t,1))=\int_{t}^1d\mu(s)$. Consider the Hankel matrix $\cH_\mu$, induced by the moment sequence 
\[
\mu_n=\int_0^1  t^n \,d\mu(t)\,.
\]
For $f(z)=\sum a_n z^n \in H(\D)$, we have 
\[
\cH_\mu(f)(z)=\sum_{n=0}^\infty \left(\sum_{k=0}^\infty \mu_{n+k}\, a_k \right) z^n\,, 
\]
and when $\mu_n=O(\frac{1}{n+1})$, i.e. $\mu(t,1)=O(1-t)$ as $t\to 1$, then 
\begin{equation}\label{integral form mu}
\cH_\mu(f)(z)=\int_0^1 f(t) \dfrac{1}{1-tz}\,d\mu(t)\,.
\end{equation}
In \cite{Galanopoulos2010}, P. Galanopoulos and J. A. Pel\'aez proved that, under the above condition on $\mu$, $\cH_\mu: H^1 \to H^1$ is bounded if and only if
\[
\mu(t,1)\,\log\dfrac{1}{1-t}=O(1-t) \text{ as }t\to 1\,.
\]
On the Bergman space, $\cH_\mu: A^2 \to A^2$ is bounded if and only if the measure $\abs{h'_\mu(z)}^2dA(z)$ is a {\it Dirichlet Carleson} measure, where $h_\mu(z)=\sum_n \mu_n\,z^n$. We note that a measure $\mu$ is Dirichlet Carleson, if the identity operator is bounded when acting from the Dirichlet space $\mathcal{D}$ to $L^2(\D, \mu)$, see \cite{Arcozzi2002}.   A few years later, C. Chatzifountas, D. Girela and Pel\'aez \cite{Chatzifountas2014} described the measures $\mu$ for which $\cH_\mu$ is bounded between different Hardy spaces, in terms of Carleson measure conditions. In 2018, Girela and N. Merch\'an \cite{Girela2018} considered the action of $\cH_\mu$ on BMOA, on the classical Bloch space $\cB$ and between different conformally invariant spaces. We recall that the Bloch space consists of the functions $f\in H(\D)$ for which $\sup_{\abs{z}<1} \abs{f'(z)}(1-\abs{z}^2)<\infty$. Recently, the study of $\cH_\mu$ was extended to Bloch-type spaces \cite{Li2021} and Zygmund spaces \cite{Manavi2022}. Furthermore, M. Beltr\'an-Meneu, J. Bonet, and E. Jord\'a in \cite{beltran2024} investigated the boundedness and compactness of $\cH_\mu$ on the Korenblum-type weighted spaces of holomorphic functions 
\[
H_v^\infty=\{f\in H(\bD):\;\sup_{z\in\bD}v(z)|f(z)|<\infty\}\,,
\]
and their "little-o" analogues $H_v^0$.
\vspace{11 pt}

A different generalization occurs if we replace the “kernel" $\frac{1}{1-tz}$ in \eqref{Integral form 1} by a more general kernel, of the form $g'(tz)$. We call $g$ the {\it symbol} of the operator $\cH_g$, and we write
\[
\cH_g(f)(z)=\int_0^1 f(t) g'(tz)\,dt\,.
\]
We notice that for the choice $g(z)=\log\frac{1}{1-z}$, $\cH_g$ reduces to $\cH$. This class of operators was studied in \cite{Galanopoulos2014} on Hardy and Bergman spaces. In specific, Galanopoulos, Girela, Pel\'aez, and Siskakis proved that $\cH_g: H^p \to H^p$, $1<p\leq 2$, is bounded if and only if the symbol $g$ is in the {\it mean Lipschitz spaces} $\Lambda(p,\frac{1}{p})$. In the case $2<p<\infty$, they proved that if $\cH_g$ is bounded, then $g \in \Lambda(p,\frac{1}{p})$ and they provided a different sufficient condition, i.e. if $g\in \Lambda(q,\frac{1}{q})$, with $q<p$, then $\cH_g$ is bounded on $H^p$. The sufficient condition is stronger than the necessary condition and they conjectured that $\cH_g$ is bounded on $H^p$ if and only if $g\in \Lambda(p,\frac{1}{p})$, for all $p>1$.

For the weighted Bergman spaces, Galanopoulos, Girela, Pel\'aez, and Siskakis proved \cite[Theorem 4]{Galanopoulos2014} that for all $p>1$ and $-1<\alpha < p-2$, the operator $\cH_g$ is bounded on $A^p_\alpha$, if and only if $g\in \Lambda(p,\frac{1}{p})$. We remind the readers that the mean Lipschitz spaces can be described via the growth of the  integral means of the derivative, that is
\[
\Lambda(p,a):=\left\{f\in H^p:\; M_p(r,f')=O\left((1-r)^{a-1}\right) \right\}\,,
\]
where $1\leq p<\infty$, and $a\in (0, 1]$.
\vspace{11 pt}

A new generalization was recently introduced in \cite{Bellavita2024} (see also \cite{Athanasiou2023}). Consider the lower triangular Hausdorff matrix $C_\mu$, induced by the moment sequence $\{\mu_n\}$ of a finite positive Borel measure $\mu$ on $(0,1)$. Its entries are given by
\[
c_{n, k}=\binom{n}{k}\int_0^1t^k(1-t)^{n-k}\,d\mu(t),
\qquad 0\leq k\leq n\,.
\]
If we eliminate  the zeros in each column of  $C_\mu$ and
shift up the columns to their  first nonzero entry, we obtain a new matrix, which is denoted by $\Gamma_\mu$, with entries given by
\[
\gamma_{n, k}=c_{n+k, k}=\binom{n+k}{k}\int_0^1t^k(1-t)^n\,d\mu(t)\,.
\]
 We note that if we apply the operation described above to the classical Ces\'aro matrix, we obtain $\cH$.
$\Gamma_\mu$ admits an integral representation given by
\[
\Gamma_\mu (f)(z)=\int_0^1 T_t f(z) \,d\mu(t)\,,
\]
where $T_t$ are the weighted composition operators in \eqref{Tt}. The authors in \cite{Bellavita2024} studied the operator $\Gamma_\mu$ on the Hardy spaces $H^p$, and proved that it induces a bounded operator, if and only if the measure $\mu$ satisfies
$$\int_{0}^{1}\frac{t^{\frac{1}{p}-1}}{(1-t)^{\frac{1}{p}}}d\mu(t),\;p>1\,,$$
and a similar condition in the case $p=1$. For $p\geq 2$, they obtained the exact value of the norm
\[
\|\Gamma_\mu\|_{H^p\to H^p}=\int_{0}^{1}\dfrac{t^{\frac{1}{p}-1}}{(1-t)^{\frac{1}{p}}}d\mu(t)\, .
\]
\noindent 
The next natural step is to consider the action of $\Gamma_\mu$ in Bergman spaces $A^p$ when $1<p<\infty$. In \cite{bellavita2024a}, Miihkinen, Norrbo, J. Virtanen, and the authors of this survey established that
\[
\|\Gamma_\mu\|_{A^p\to A^p}=\int_{0}^{1}\dfrac{t^{\frac{2}{p}-1}}{(1-t)^{\frac{2}{p}}}\,d\mu(t)\,,\;\;4\leq p<\infty\,,
\]
 and they also characterized the measures $\mu$ for which the operator is bounded in the range $1\leq p<4$. 
%In an upcoming article, we show that, using techniques similar to the Hardy case, one can compute the exact norm of $\Gamma_\mu$ for $p\geq 4$, by considering some natural integrability condition on the measure $\mu$.More precisely\[\|\Gamma_\mu\|_{A^p\to A^p}=\int_{0}^{1}\dfrac{t^{\frac{2}{p}-1}}{(1-t)^{\frac{2}{p}}}\,d\mu(t)\,,\;\;4\leq p<\infty\,,\]under the assumption that the above integral is finite.
%%%%%%%%%%%%%%%%%%%%%%%%%%%%%%%%%%%%%%%%%%%%%%%%%%%%%%%%%%%%%%
%\section{A multiplicative Hankel analogue}\noindent
\vspace{11 pt}

An analogue of a Hankel matrix defined on the infinite-dimensional torus $\mathbb{T}^\infty$,  is an infinite matrix whose entries $a_{n,k}$ depend only on the product $n\cdot k$. We call this Hankel forms Helson matrices, after H. Helson (see \cite{Helson2005,Helson2006,Helson2010}) who studied the problem of
whether every multiplicative Hankel form on $\mathbb{T}^\infty$, comes from a bounded function.

In 2016, O. F. Brevig, K.-M. Perfekt, K. Seip, Siskakis and Vukoti\'c \cite{Brevig2016} considered a multiplicative Hilbert matrix $\mathbf{M}$, given by
\[\label{eq:matrixmultiplicative}
\mathbf{M} \,=\,
\left(\frac{1}{\sqrt{nk}\log(nk)}\right)_{n,k\,\geq\, 2}\, .
\]
They showed that $\mathbf{M}$ admits an integral representation on the space 
$\ell^2(\mathbb{N}_{\geq 2})$. In specific, if $f(s) =\sum_{n\geq 2} c_n\, n^{-s}$ with $\Re(s)>1/2$, then
\[
\mathbf{M}(f)(s) =\sum_{n\geq 2}\left(\sum_{k\geq 2} \mathbf{M}_{n,k}\, c_k \right)\, \frac{1}{n^s}= \int_{1/2}^{\infty} f(w)\left( \zeta(w+s)-1\right)\, dw\,,
\]
where $\zeta(s)$ denotes the classical Riemann's zeta function. The function $f$ defined in this way belongs to the Hilbert space of Dirichlet series $\mathcal{H}^2_0$. We refer the interested readers to \cite{Hedenmalm1997} for an introduction to Hardy spaces of Dirichlet series. They were able to prove that the operator $\mathbf{M}$ is bounded on $\mathcal{H}^2_0$ with operator norm equal to $\pi$. For $1<p<\infty$, they also considered 
\[
\mathbf{M}_p=\left( \frac{1}{n^{\frac{p-1}{p}}k^{\frac{1}{p}}\log(nk)}\right)_{n,k\, \geq\, 2}\ 
\]
acting on $\ell^p(\mathbb{N}_{\geq 2})$, obtaining that $\|M_p\|=   \pi/\sin(\pi/ p)$. 

The analogy with the classical Hilbert matrix
breaks down when one considers $\mathcal{H}^q_0$ with $q\neq 2$. Indeed, they verified that $\mathbf{M}$ cannot be extended to be a bounded operator on the $H^q$ analogues of $\mathcal{H}^2_0$, which by Bayart’s work \cite{Bayart2002} are associated with $H^q(\mathbb T^\infty)$. This negative result is related to the fact
that $H^q(\mathbb T^\infty)$ is not complemented in $L^q(\mathbb T^\infty)$.

\section{Spectral properties}
\noindent In 1950, W. Magnus \cite{Magnus1950} studied the spectrum of the Hilbert matrix acting on $\ell^2$: he proved that the spectrum of $\cH$ is purely continuous and every real value of $\gamma$ with $0\leq \gamma \leq \pi$ belongs to its spectrum. 

In the following years, a slightly more general version of 
$\mathbb{H}$ has been studied, namely $\cH_\lambda$,  with entries given by $\frac{1}{n+k+\lambda}$ with $\lambda \in \mathbb{C}$. The research concentrated on finding {\it latent roots}, i.e. “eigenvalues" whose corresponding eigenvector does not belong to the considered space. 
In specific, T. Kato \cite{Kato1957} proved that $\pi$ is a latent root of $\cH$ acting in $\ell^2$. 
Expanding the work of Magnus, M. Rosenblum \cite{Rosenblum1958,Rosenblum1958b} showed that when $\lambda\in \mathbb{R}\setminus \mathbb{Z}$, then $\cH_\lambda: \ell^2 \to \ell^2$ has a continuous spectrum of multiplicity $1$ inside $[0,\pi]$. Furthermore, if $\lambda\geq 1/2$, there are no eigenvalues, while, if $0<\lambda<1/2$, then $\cH_\lambda$ has only two eigenvalues: $\frac{\pm \pi}{\sin[\pi(1-\lambda)]}$. 
Moreover, C. K. Hill \cite{Hill1960,Hill1961} proved that for each complex number $\mu$ with $0<Re(\mu)<1/2$, the quantity $\pi/\sin(\pi\mu)$ is a latent root of $\cH_\lambda$.

For a number of years no significant progress was made, until 2012, when A. Aleman, A. Montes-Rodriguez, and A. Sarafoleanu \cite{Aleman2012}, inspired by the work of Rosenblum and Hill, studied in depth the spectrum of 
\[
\cH_\lambda=\left( \frac{1}{n+k+\lambda}\right),\;\lambda \in \mathbb{C}\setminus \mathbb{Z}\,.
\]
They showed that in $H^p,\;1<p<2$, and for each $\alpha$ with $\frac{1}{2}<Re(\alpha)<\frac{1}{p}$, we have that $\frac{\pi}{\sin\pi\alpha}$ is an eigenvalue for $\cH_\lambda$, with corresponding eigenfunction
\[
g_\alpha(z)=(1-z)^\alpha\; _2F_1(\alpha+1, \alpha+\lambda;\lambda;z)\,.
\]
Their approach was based on certain differential operators that “almost" commute with $\cH_\lambda$ and on properties of hypergeometric functions.

A major breakthrough was made by B. Silbermann's work \cite{Silbermann2021}, who was able to describe the spectrum of $\cH$, when acting on the Hardy spaces $H^p$, and the sequence spaces $\ell^p$. In specific, he proves the following results \cite[Theorem 6.1]{Silbermann2021}:

\noindent 
For $\lambda \in \mathbb{C}\setminus \mathbb{Z}_{\leq 0}$ and $\cH_\lambda: H^p \to H^p,\;1<p<\infty$,
\begin{enumerate}
    \item The spectrum of $\cH_\lambda$ is equal to the image of the function 
    $$
  V_p(\xi)=i\pi \sinh^{-1}\left(\pi(\xi+\frac{i}{p})\right),\;\xi\in \mathbb{R}\,.
    $$
    \item If $N$ is the largest integer for which $(1-z)^{-\lambda - N}\in H^p,\; p\geq 2$, then the functions
    \begin{equation}\label{eigenfunction}
    f_n(z)=(1-z)^{-\lambda - n}(1+z)^n,\;0\le n \le N\,,
\end{equation}
    are eigenfunctions of $\cH_\lambda$, with corresponding eigenvalues 
    $$(-1)^n \dfrac{\pi}{\sin(\lambda \pi)}\,.$$
    \item $\cH_\lambda$ may have a one-dimensional eigenspace, with eigenvalues $\frac{\pm \pi}{\sin(\lambda \pi)}$, which is not contained in $H^2$.
\end{enumerate}
For the $\ell^p$ setting, exploiting the Hausdorff-Young inequalities, he obtained analogous results \cite[Theorem 7.6]{Silbermann2021}. That is, the spectrum of $\cH_\lambda: \ell^p \to \ell^p$, $1<p<\infty$,  is equal to the image of the function $ V_q(\xi)$. For $q\geq2$, with $\frac{1}{p}+\frac{1}{q}=1$, the Fourier coefficients of the functions $f_n$ defined in \eqref{eigenfunction}, are eigenfunctions of $\cH_\lambda$ acting on $\ell^q$, related to the same eigenvalues, as in the Hardy case.
His approach is based on “heavy" mathematical tools, such as the theory of Fredholm operators, the Gelfand transform, the theory of Toeplitz and Hankel operators, and a variety of other techniques, thus leading to an extremely non-trivial proof of the results.

Recently, Montes-Rodriguez and Virtanen in \cite{Montes2024} were able to provide simpler proofs of the classical spectral results of Magnus and Hill for the Hilbert matrix operator $\mathcal{H}$, utilizing the Mehler-Fock transform, that is,
\[
f \mapsto \mathcal{P}(f)(t):=\int_1^{\infty}f(x)P_{it-1/2}(x)dx\ , t\geq 0\, ,
\]
where $P_\mu$ is a Legendre function, see \cite{Lebedev1965}. They also showed that the measure
\[
dr(t)=\frac{2}{\pi^2}\text{arcosh}\left(\frac{\pi}{x}\right)dx\, \text{ with } x \in [0,\pi]
\]
is the spectral measure of $\mathcal{H}$.

On the setting of the multiplicative Hilbert matrix, Brevig et al. \cite{Brevig2016} showed that $\mathbf{H}$ has no eigenvalues and its continuous spectrum is $[0, \pi]$. Furthermore, Perfekt and A. Pushnitski \cite{Perfekt2018} proved  that the multiplicative Hilbert matrix has no singular continuous spectrum and that its absolutely continuous spectrum has multiplicity one. More precisely, they studied a more general class of Helson matrices, corresponding to the sequence
\[
g_a(n)=\dfrac{1}{\sqrt{n}(a+\log n)},\;n\geq 1,\;a\geq 0\,.
\]
For $a=0$, the matrix reduces to $\mathbf{H}$.
%%%%%%%%%%%%%%%%%%%%%%%%%%%%%%%%%%%%%%%%%%%%%%%%%%%%%%%%%%%%%%%%%%%%%%%%%%%%%%%%%%%%%%%%%%%%%

As a final remark for this section, we would like to discuss the anticipated upcoming article by Aleman, Siskakis, and Vukoti\'c \cite{Aleman2024}, that is expected to shed light on the spectrum of the classical Hilbert matrix on a broad range of holomorphic function spaces. 
This class consists of the weighted conformally invariant spaces of analytic functions on $\D$, i.e. for each such space $X$ the weighted composition operator
\[
W_\phi(f)(z)=\left(\phi'(z)\right)^\alpha f(\phi(z)),\; 0<\alpha<1
\]
is bounded on $X$ and $\sup_{\phi}\|W_\phi\|_{X\to X}<\infty$, where $\phi$ is a M\"obius authomorphism of the disk. This class includes Hardy, Bergman and weighted Dirichlet spaces of appropriate parameter. 

Using the results obtained on  spaces satisfying these requirements, they are able to deduce information for other spaces, including the classical $\ell^p$ sequence spaces. 
Their approach is based on the study of a Carleman-type \textit{companion operator}, defined by
\[
\mathcal{K}(f)(z)=\int_{-1}^1 f(t)\dfrac{1}{1-tz}\,dt\,, 
\]
with a corresponding matrix which is a “reduced" Hilbert matrix, whose entries are zero when $n+k+1$ is an even number. They are able to identify the spectrum of $\mathcal{K}$ and its various parts, and this allows to conclude analogous results for $\cH$. 

%%%%%%%%%%%%%%%%%%%%%%%%%%%%%%%%%%%%%%%%%%%%%%%%%%%%%%%%%%%%%%
\section{Acknowledgments}
\noindent The authors would like to thank professor A. G. Siskakis for introducing them to this fascinating topic, and for all his help during the preparation of this survey.

%%%%%%%%%%%%%%%%%%%%

\bibliography{Literature}
\bibliographystyle{plain}

\end{document}